\documentclass[reqno,11pt]{amsart}
\usepackage[english]{babel}
\usepackage{amssymb,amsthm,amsfonts}
\usepackage[T1]{fontenc}
\usepackage[latin1]{inputenc}
\usepackage[dvips]{graphicx}
\usepackage{fancyhdr}
\usepackage{color}
\usepackage{textcomp}

\allowdisplaybreaks[3]
\renewcommand{\div}{\mathop{\rm div}\nolimits}
\renewcommand{\det}{\mathop{\rm det}\nolimits}
\newcommand{\tr}{\mathop{\rm tr}\nolimits}

\newcommand{\Sym}{\mathcal{S}}
\newcommand{\R}{\mathbb{R}}
\newcommand{\N}{\mathbb{N}}
\newcommand{\C}{\mathbb{C}}
\newcommand{\Z}{\mathbb{Z}}

\newcommand{\dt}{\partial_t}
\providecommand{\abs}[1]{\left\lvert#1\right\rvert}
\providecommand{\norm}[1]{\left\lVert#1\right\rVert}

\newenvironment{dimo}{\begin{proof}[\textbf{Proof}]}{\end{proof}}
\newenvironment{claim}{\begin{proof}[{CLAIM}]}{\end{proof}}

\newtheorem{teor}{Theorem}[section]
\newtheorem{lemm}[teor]{Lemma}
\newtheorem{defin}[teor]{Definition}
\newtheorem{propos}[teor]{Proposition}

\newtheorem{obs}[teor]{Remark}

\begin{document}

\author{Elisabetta Chiodaroli}
\title[A counterexample to well-posedness
 of entropy solutions]{A counterexample to well-posedness \\
 of entropy solutions \\
to the compressible Euler system}

\begin{abstract}
We deal with entropy solutions to the Cauchy problem
for the isentropic compressible Euler equations in the space-periodic case.
In more than one space dimension, the methods developed by De Lellis-Sz\'{e}kelyhidi enable us to show failure of uniqueness on a finite time-interval
for entropy solutions starting from any continuously differentiable initial 
density and suitably constructed bounded initial linear momenta. 
\end{abstract}

\maketitle 

\pagestyle{myheadings}

\section{Introduction}
\setcounter{equation}{0}

In this note, we deal with the Cauchy Problem for the isentropic compressible
Euler equations in the space-periodic setting. Given any continuously differentiable initial density, we can
construct bounded initial linear momenta
for which admissible solutions are not
unique in more than one space dimension. 

We first introduce the isentropic compressible Euler
equations of gas dynamics in $n$ space dimensions, $n\geq 2$ (cf. Section 3.3 of \cite{da}). They are obtained as a simplification
of the full compressible Euler equations, by assuming the entropy to be constant.
The state of the gas will be described through the state vector
$$ V=(\rho, m)$$
whose components are the density $\rho$ and the linear momentum $m$. The balance laws in force are for mass
and linear
momentum. The resulting system, which consists of $n+1$ equations,
takes the form:
\begin{equation}\label{eq:isentropic Euler system}
\left\{\begin{array}{l}
\partial_t \rho + {\rm div}_x m \;=\; 0\\
\partial_t m + {\rm div}_x \left(\frac{m\otimes m }{\rho}\right) + \nabla_x [ p(\rho)]\;=\; 0\\
\rho (\cdot,0)\;=\; \rho^0\\
m (\cdot, 0)\;=\; m^0 \,
\end{array}\right..
\end{equation}
The pressure $p$ is
a function of $\rho$ determined from the
constitutive thermodynamic relations of the gas in question. A
common choice is the polytropic pressure law
$$p(\rho)= k \rho^\gamma$$
with constants $k>0$ and $\gamma>1$. The set of admissible values is $P=\{\rho>0\}$ (cf. \cite{da} and \cite{se}).
The system is hyperbolic if
$$ p' (\rho)>0.$$
In addition, thermodynamically admissible processes must also satisfy
an additional constraint coming from the energy inequality
\begin{equation} \label{eq:energy inequality}
\dt \left(\rho \varepsilon(\rho)+\frac{1}{2}
\frac{\abs{m}^2}{\rho}\right)+\div_x
\left[\left(\varepsilon(\rho)+\frac{1}{2}
\frac{\abs{m}^2}{\rho^2}+\frac{p(\rho)}{\rho}\right) m \right]
\;\leq\; 0
\end{equation}
where the internal energy $\varepsilon: \R^+\rightarrow \R$ is given
through the law $p(r)=r^2 \varepsilon'(r)$. The physical region for
\eqref{eq:isentropic Euler system} is $\{(\rho, m) | \abs{m}\leq R
\rho\}$, for some constant $R>0$. For $\rho>0$, $v=m/\rho$ represents
the velocity of the fluid.


We will consider, from now on, the case
of general pressure laws given by a function $p$ on $[0,\infty[$,
that we always assume to be continuously differentiable on $[0,\infty[$.
The crucial requirement we
impose upon $p$ is that it has to be strictly increasing on $[0,\infty[$.
Such a condition is meaningful from a physical viewpoint since it is
a consequence of the principles of thermodynamics.

Now, we will disclose the content of this note.
Using some techniques introduced by De Lellis-Sz\'{e}kelyhidi (cf. \cite{dls1} and \cite{dls2}) we can consider
any continuously differentiable periodic initial density $\rho^0$ and exhibit suitable periodic initial momenta $m^0$ for which space-periodic weak
admissible solutions of \eqref{eq:isentropic Euler system} are not
unique on some finite time-interval.

\begin{teor} \label{th:isentropic non-uniqueness} 
Let $n\geq 2$. Then, for any given function $p$ and any given continuously differentiable periodic initial density $\rho^0$,
there exist a bounded periodic initial momentum $m^0$ and a positive time $\overline{T}$ for which there are infinitely many space-periodic admissible solutions $(\rho,m)$
of \eqref{eq:isentropic Euler system} on $\R^n\times [0,\overline{T}[$ with $\rho\in C^1(\R^n\times [0,\overline{T}[)$.
\end{teor}

\begin{obs}
Indeed, in order to prove Theorem \ref{th:isentropic non-uniqueness}, it would be enough to assume that the initial density is a H\"older continuous periodic function:
$\rho^0\in C^{0,\alpha}(\R^n)$ (cf. Proof of Proposition \ref{prop:starting subsolution}).
\end{obs}

Some connected results are obtained in \cite{dls2} (cf. Theorem 2 therein) as a further consequence of their 
analysis on the incompressible Euler equations. Inspired by their approach, we adapt and apply directly to \eqref{eq:isentropic Euler system}
the method of convex integration combined with Tartar's programme on oscillation phenomena in conservation laws (see \cite{ta} and \cite{kms}). In this way, we can show failure of uniqueness
of admissible solutions to the compressible Euler equations starting from any given continuously differentiable initial density. For a survey on these \textit{h}-principle-type results
in fluid dynamics we refer the reader to \cite{dls3}.

The paper is organised as follows. Section 2 is an overview on the definitions of weak and admissible solutions to \eqref{eq:isentropic Euler system} and 
gives a first glimpse on how Theorem \ref{th:isentropic non-uniqueness} is achieved. Section 3 is devoted to the reformulation of a simplified version of the isentropic compressible Euler
equations as a differential inclusion and to the corresponding geometrical analysis.
In Section 4 we state and prove a criterion (Proposition \ref{prop:criterion}) to select initial momenta allowing for infinitely many solutions. The proof builds upon a refined version of the Baire
category method for differential inclusions developed in \cite{dls2} and aimed at yielding weakly continuous in time solutions. Section 5 and 6 contain the proofs of the main 
tools used to prove Proposition \ref{prop:criterion}. In Section 7, we show initial momenta satisfying the requirements of Proposition \ref{prop:criterion}.
Finally, in Section 8 we prove Theorem \ref{th:isentropic non-uniqueness} by applying Proposition \ref{prop:criterion}.

\section{Weak and admissible solutions to \\ the isentropic Euler system}
\setcounter{equation}{0}

The deceivingly simple-looking system of first-order partial differential
equations \eqref{eq:isentropic Euler system} has a long history of important
contributions over more than two centuries. We recall a few classical facts on this system (see for instance \cite{da} for more details).
\begin{itemize}
\item If $\rho^0$ and $m^0$ are ``smooth'' enough (see Theorem 5.3.1 in \cite{da}),
      there exists a maximal time interval $[0,T[$ on which there exists a unique ``smooth''
      solution $(\rho, m)$ of \eqref{eq:isentropic Euler system} (for $0\leq t<T$). In addition, if
      $T<\infty$, and this is the case in general, $(\rho, m)$ becomes discontinuous as $t$ goes to $T$.

\item If we allow for discontinuous solutions, i.e., for instance, solutions
      $(\rho, m) \in L^\infty$ satisfying \eqref{eq:isentropic Euler system}
      in the sense of distributions,
      then solutions are neither unique nor stable. More precisely, one can exhibit
      sequences of such solutions which  converge weakly in $L^\infty-\ast$ to functions
      which do not satisfy \eqref{eq:isentropic Euler system}.

\item In order to restore the stability of solutions and (possibly) the uniqueness,
      one may and should impose further restrictions on bounded solutions of \eqref{eq:isentropic Euler system},
      restrictions which are known as (Lax) entropy inequalities.
\end{itemize}

This note stems from the problem of better understanding the efficiency of
entropy inequalities as selection criteria among weak solutions.

Here, we have chosen to emphasize the case of the flow with \textit{space periodic} boundary conditions.
For \textit{space periodic} flows we assume that the fluid fills the entire space $\R^n$ but with the condition that $m,\rho$ are periodic functions 
of the space variable. The space periodic case is not a physically achievable one, but it is relevant on the physical side as a model for some flows.
On the mathematical side, it retains the complexities due to the nonlinear terms (introduced by the kinematics) and therefore it includes many
of the difficulties encountered in the general case.
However the former is simpler to treat because of the absence of boundaries. Furthermore, using Fourier transform as a tool simplifies the analysis.

Let $Q=[0,1]^n$, $n\geq 2$ be the unit cube in $\R^n$. We denote by $H^m_p(Q)$, $m\in \N$, the space of functions which are
in $H^m_{loc}(\R^n)$ and which are periodic with period $Q$:
\begin{equation*}
 m(x+l)=m(x) \quad \text{for a.e. } x\in\R^n \text{ and every } l\in \Z^n.
\end{equation*}
For $m=0$, $H^0_p(Q)$ coincides simply with $L^2(Q)$. Analogously, for every functional space $X$ we define $X_p(Q)$ to be the space of functions 
which are locally (over $\R^n$) in $X$
and are periodic of period $Q$. 
The functions in $H^m_p(Q)$ are easily characterized by their Fourier series expansion
\begin{equation} \label{eq:Fourier expansion}
 H^m_p(Q)=\left\{m\in L^2_p(Q): \sum_{k\in\Z^n} \abs{k}^{2m} \abs{\widehat{m}(k)}^2< \infty \text{ and } \widehat{m}(0)=0\right\},
\end{equation}
where $\widehat{m}: \Z^n\rightarrow \C^n$ denotes the Fourier transform of $m$. We will use the notation $H(Q)$ for $H^0_p(Q)$ and $H_w(Q)$ for the space $H(Q)$
endowed with the weak $L^2$ topology.

Let $T$ be a fixed positive time. By a \textit{weak solution} of \eqref{eq:isentropic Euler system} on $\R^n\times[0,T[$ we
mean a pair $(\rho, m)\in L^\infty([0,T[ ;L^\infty_p(Q))$ satisfying
\begin{equation} \label{eq:bound on m}
 \abs{m(x,t)}\leq R \rho(x,t) \qquad \text{for a.e. } (x,t)\in \R^n\times [0,T[ \text{ and some } R>0,
\end{equation}
and such that the following identities hold for every test functions $\psi\in C_c^{\infty}([0, T[ ; C^{\infty}_p (Q))$,
$\phi\in C_c^{\infty}([0, T[ ; C^{\infty}_p (Q))$:
\begin{equation} \label{eq:weak1}
\int_0^{T} \int_{Q} \left[\rho\dt \psi+ m \cdot \nabla_x \psi\right] dx dt+\int_{Q} \rho^0(x)\psi(x,0) dx \;=\; 0
\end{equation}
\begin{align} \label{eq:weak2}
&\int_0^{T} \int_{Q} \left[ m \cdot \dt \phi+ \left\langle \frac{m \otimes m}{\rho}, \nabla_x \phi\right\rangle +p(\rho) \div_x \phi \right] dx dt\notag\\
&+\int_{Q} m^0(x)\cdot\phi(x,0) dx\;=\; 0.
\end{align}

For $n\geq 2$ the only non-trivial entropy is the total energy $\eta=\rho \varepsilon(\rho)+\frac{1}{2}
\frac{\abs{m}^2}{\rho}$ which corresponds to the flux $\Psi=\left(\varepsilon(\rho)+\frac{1}{2}
\frac{\abs{m}^2}{\rho^2}+\frac{p(\rho)}{\rho}\right) m$.

Then a bounded weak solution $(\rho, m)$ of \eqref{eq:isentropic Euler system} satisfying \eqref{eq:energy inequality}
in the sense of distributions, i.e. satisfying the following inequality
\begin{align} \label{eq:admissibility condition}
 &\int_0^{T}\int_{Q} \left[\left(\rho\varepsilon(\rho)+\frac{1}{2}
\frac{\abs{m}^2}{\rho}\right)\dt \varphi+\left(\varepsilon(\rho)+\frac{1}{2}
\frac{\abs{m}^2}{\rho^2}+\frac{p(\rho)}{\rho}\right) m \cdot \nabla_x \varphi \right]\notag \\
&+\int_{Q} \left(\rho^0\varepsilon(\rho^0)+\frac{1}{2}
\frac{\abs{m^0}^2}{\rho}\right) \varphi(\cdot,0) \;\geq\; 0,
\end{align}
for every nonnegative $\varphi\in C_c^{\infty}([0,T[ ; C^{\infty}_p (Q))$, is said to be an \textit{entropy} (or \textit{admissible}) solution of \eqref{eq:isentropic Euler system}.

The lack of entropies is one of the essential reasons for a very limited understanding of compressible
Euler equations in dimensions greater than or equal to $2$.

A recent paper by De Lellis-Sz\'{e}kelyhidi gives an example in
favour of the conjecture that entropy solutions to the
multi-dimensional compressible Euler equations are in general not
unique. Showing that this conjecture is true has far-reaching
consequences. The entropy condition is not sufficient as a
selection principle for physical/unique solutions. The
non-uniqueness result by De Lellis-Sz\'{e}kelyhidi is a byproduct of
their new analysis of the incompressible Euler equations based on
its formulation as a differential inclusion. They first show that,
for some bounded compactly supported initial data, none of the
classical admissibility criteria singles out a unique solution to
the Cauchy problem for the incompressible Euler equations. As a
consequence, by constructing a piecewise
constant in space and independent of time density $\rho$, they look at the
compressible isentropic system as a ``piecewise incompressible''
system (i.e. still incompressible in the support of the velocity field) and thereby exploit the result for the incompressible Euler
equations to exhibit bounded initial density and bounded compactly supported initial momenta for which admissible
solutions of \eqref{eq:isentropic Euler system} are not unique (in
more than one space dimension).

Inspired by their techniques, we give a further counterexample to the well-posedeness of
entropy solutions to \eqref{eq:isentropic Euler system}. Our result differs in two main aspects: here the initial density can be any given ``regular'' function and remains ``regular''
forward in time while in \cite{dls2} the density allowing for infinitely many admissible solutions must be chosen as piecewise constant in space; on the other hand we are 
not able to deal with compactly supported momenta (indeed we work in the periodic setting), hence our non-unique entropy solutions are only locally $L^2$ in contrast with the global-$L^2$-in-space property of solutions obtained
in \cite{dls2}. Moreover, we have chosen to study the case of the flow in a cube of $\R^n$ with \textit{space periodic} boundary conditions.
This case leads to many technical simplifications while retaining the main structure of the problem.

More precisely, we are able to analyze the compressible Euler equations in the framework of convex integration. This method works well 
with systems of nonlinear PDEs such that the convex envelope (in an appropriate sense) of each small domain of the submanifold representing the PDE in the jet-space (see \cite{em} for more details)
is big enough. In our case, we consider a simplification of system \eqref{eq:isentropic Euler system}, namely the semi-stationary associated problem, whose submanifold
allows a convex integration approach leading us to recover the result of Theorem \ref{th:isentropic non-uniqueness}.

We are interested in the semi-stationary Cauchy problem associated
with the isentropic Euler equations (simply set to $0$ the time
derivative of the density in \eqref{eq:isentropic Euler system} and
drop the initial condition for $\rho$):
\begin{equation}\label{eq:semistationary system}
\left\{\begin{array}{l}
{\rm div}_x m \;=\; 0\\
\partial_t m + {\rm div}_x \left(\frac{m\otimes m }{\rho}\right) + \nabla_x [ p(\rho)]\;=\; 0\\
m (\cdot, 0)\;=\; m^0. \
\end{array}\right.
\end{equation}
A pair $(\rho, m)\in L^\infty_p(Q)\times L^{\infty}([0, T[; L^\infty_p(Q))$ is a \textit{weak solution} on $\R^n\times [0,T[$ of
\eqref{eq:semistationary system} if $m(\cdot, t)$ is weakly-divergence free for almost every $0<t<T$ and satisfies the following bound
\begin{equation} \label{eq:bound on m2}
 \abs{m(x,t)}\leq R \rho(x) \qquad \text{for a.e. } (x,t)\in \R^n\times [0,T[ \text{ and some } R>0,
\end{equation}
and if the following identity holds for every $\phi\in C_c^{\infty}([0, T[ ; C^{\infty}_p (Q))$:
\begin{align} \label{eq:semistationaryweak2}
&\int_0^{T}\int_{Q} \left[m \cdot \dt \phi+ \left\langle\frac{m\otimes m}{\rho}, \nabla_x \phi\right\rangle +p(\rho) \div_x \phi\right] dx dt \notag\\
&+\int_{Q} m^0(x)\cdot\phi(x,0) dx\;=\; 0.
\end{align}

A general observation suggests us that a non-uniqueness result for weak solutions of \eqref{eq:semistationary system} whose momentum's magnitude satisfies some
suitable constraint could lead us to a non-uniqueness result for entropy solutions of the isentropic Euler equations \eqref{eq:isentropic Euler system}.
Indeed, the entropy solutions we construct in Theorem \ref{th:isentropic non-uniqueness}
come from \textit{some} weak solutions of
\eqref{eq:semistationary system}.

\begin{teor} \label{th:semistationary non-uniqueness}
Let $n\geq 2$. Then, for any given function $p$, any given density $\rho_0\in C^1_p(Q)$ and any given finite positive time $T$,
there exists a bounded initial momentum $m^0$ for which there are infinitely many weak solutions $(\rho,m)\in C^1_p(Q)\times C([0,T];H_w(Q))$
of \eqref{eq:semistationary system} on $\R^n\times [0,T[$ with density $\rho(x)=\rho_0(x)$. \\
In particular, the obtained weak solutions $m$ satisfy
\begin{align} 
 &\abs{m(x,t)}^2=\rho_0(x) \chi(t) \quad \text{a.e. in }\R^n \times [0,T[, \label{eq:modulus constraint a}\\
&\abs{m^0(x)}^2=\rho_0(x) \chi(0) \quad \text{a.e. in }\R^n, \label{eq:modulus constraint b}
\end{align}
for some smooth function $\chi$.

\end{teor}

An easy computation shows how, by properly choosing the function $\chi$ in \eqref{eq:modulus constraint a}-\eqref{eq:modulus constraint b}, the solutions $(\rho_0,m)$ 
of \eqref{eq:semistationary system} obtained in Theorem \ref{th:semistationary non-uniqueness} satisfy the admissibility condition \eqref{eq:admissibility condition}.

\begin{teor} \label{th:function chi}
Under the same assumptions of Theorem \ref{th:semistationary non-uniqueness}, there exists a maximal time $\overline{T}>0$ such that
the weak solutions $(\rho,m)$ of \eqref{eq:semistationary system} (coming from Theorem \ref{th:semistationary non-uniqueness})
satisfy the admissibility condition \eqref{eq:admissibility condition} on $[0,\overline{T}[$.
\end{teor}

Our construction yields initial data $m^0$ for which the
nonuniqueness result of Theorem \ref{th:isentropic non-uniqueness}
holds on any time interval $[0,T[$, with $T\leq \overline{T}$. However, as pointed out before, for
sufficiently regular initial data, classical results give the local
uniqueness of smooth solutions. Thus, \textit{a fortiori}, the
initial momenta considered in our examples have necessarily a certain
degree of irregularity.

\section{Geometrical analysis}
\setcounter{equation}{0}
This section is devoted to a qualitative analysis of the isentropic compressible Euler equations
in a semi-stationary regime (i.e. \eqref{eq:semistationary system}).

As in \cite{dls1} we will interpret the system \eqref{eq:semistationary system}
in terms of a differential inclusion, so that it can be studied in the framework combining
the plane wave analysis of Tartar, the convex integration of Gromov and the Baire's arguments.

\subsection{Differential inclusion}
The system \eqref{eq:semistationary system} can indeed be naturally expressed as a linear system
of partial differential equations coupled with a pointwise nonlinear constraint, usually called \textit{differential inclusion}.

The following Lemma, based on Lemma $2$ in \cite{dls2}, gives such a reformulation. We will denote by $\Sym^n$ the space of symmetric $n\times n$ matrices, by $\Sym^n_0$
the subspace of $\Sym^n$ of matrices with null trace, and by $I_n$ the $n\times n$ identity matrix.

\begin{lemm} \label{lemma:differential inclusion}
Let $m\in L^{\infty}([0,T]; L^\infty_p(Q;\R^n))$, $U\in L^{\infty}([0,
T]; L^\infty_p (Q;\mathcal{S}_0^n))$ and $q\in L^{\infty}([0,T];L^\infty(Q;\R^+))$ such that
\begin{align} \label{eq:linearized system}
\div_x m & \;=\;0 \notag \\
\dt m+\div_x U+\nabla_x q &\;=\;0.
\end{align}
If $(m, U, q)$ solve \eqref{eq:linearized system} and in addition
there exists $\rho \in L^\infty_p(\R^n;\R^+)$ such that \eqref{eq:bound on m2} holds and
\begin{align} \label{eq:pointwise constraint}
U=\frac{m\otimes m}{\rho}-\frac{\abs{m}^2}{n\rho} I_n \quad
\mbox{a.e. in } \R^n\times [0,T],\notag\\
q=p(\rho)+\frac{\abs{m}^2}{n \rho} \quad \mbox{a.e. in } \R^n\times [0,T],
\end{align}
then $m$ and $\rho$ solve \eqref{eq:semistationary system} distributionally.
Conversely, if $m$ and $\rho$ are weak solutions of \eqref{eq:semistationary
system}, then $m$, $U=\frac{m\otimes
m}{\rho}-\frac{\abs{m}^2}{n\rho}I_n$ and
$q=p(\rho)+\frac{\abs{m}^2}{n \rho}$ solve \eqref{eq:linearized
system}-\eqref{eq:pointwise constraint}.
\end{lemm}

In Lemma \ref{lemma:differential inclusion} we made clear the
distinction between the augmented system \eqref{eq:linearized
system}, whose linearity allows a plane wave analysis, and the
nonlinear pointwise constraint \eqref{eq:pointwise constraint}, which leads us to study the graph
below.

For any given $\rho\in ]0, \infty[$, we define the following graph
\begin{align} \label{eq:K}
K_{\rho}:=\Bigg\{&(m, U, q)\in \R^n\times \mathcal{S}_0^n\times
\R^+: U=\frac{m\otimes m}{\rho}-\frac{\abs{m}^2}{n\rho}I_n, \notag \\
 q&=p(\rho)+\frac{\abs{m}^2}{n \rho}\Bigg\}.
\end{align}

The key of the forthcoming analysis is the behaviour of the graph $K_{\rho}$ with respect to the
wave vectors associated with the linear system \eqref{eq:linearized system}: are differential and algebraic
constraints in some sense compatible?

For our purposes, it is convenient to consider ``slices'' of the graph $K_{\rho}$, by considering vectors $m$
whose modulus is subject to some $\rho$-depending condition. Thus, for any given $\chi\in\R^+$, we define:

\begin{align} \label{eq:Kf}
K_{\rho,\chi}:=\Bigg\{&(m, U, q)\in \R^n\times \mathcal{S}_0^n\times
\R^+: U=\frac{m\otimes m}{\rho}-\frac{\abs{m}^2}{n\rho}I_n, \notag
\\  q&=p(\rho)+\frac{\abs{m}^2}{n \rho},\abs{m}^2=\rho \chi \Bigg\}.
\end{align}

\subsection{Wave cone}

Following Tartar's framework \cite{ta}, we consider a system of first order linear PDEs
\begin{equation} \label{eq:linear PDE}
 \sum_i A_i \partial_i z=0
\end{equation}
where $z$ is a vector valued function and the $A_i$ are matrices.
Then, \textit{planewave} solutions to \eqref{eq:linear PDE} are solutions of the form
\begin{equation} \label{eq:planewave solutions}
 z(x)= a h (x \cdot \xi),
\end{equation}
with $h:\R\rightarrow \R$. In order to find such solutions, we have to solve the relation $\sum_i \xi_i A_i a=0$,
where $\xi_i$ is the oscillation frequency in the direction $i$. The set of directions $a$ for which a solution $\xi\neq 0$
exists is called \textit{wave cone} $\Lambda$ of the system \eqref{eq:linear PDE}: equivalently $\Lambda$ characterizes the directions of one dimensional
high frequency oscillations compatible with \eqref{eq:linear PDE}.

The system \eqref{eq:linearized system} can be analyzed in this
framework. Consider the $(n+1)\times(n+1)$ symmetric matrix in block
form
\begin{equation} \label{eq:matrix}
    M=\left( \begin{array}{cc}
    U+q I_n & m \\
    m & 0\\
    \end{array} \right).
\end{equation}
Note that, with the new coordinates $y=(x,t)\in \R^{n+1}$, the system \eqref{eq:linearized system}
can be easily rewritten as $\div_y M=0$ (the divergence of $M$ in space-time is zero). Thus,
the wave cone associated with the system \eqref{eq:linearized system} is equal to
\begin{equation} \label{eq:wave cone}
 \Lambda=\left\{(m, U, q)\in \R^n\times\Sym_0^n\times \R^+: \det \left( \begin{array}{cc}
    U+q I_n & m \\
    m & 0\\
    \end{array} \right)=0\right\}.
\end{equation}

Indeed, the relation $\sum_i \xi_i A_i a=0$ for the system \eqref{eq:linearized system}
reads simply as $M\cdot (\xi, c)=0$, where $(\xi, c)\in \R^n \times \R$ ($\xi$ is the space-frequency and $c$
the time-frequency): this equation admits a non-trivial solution if $M$ has null determinant,
hence \eqref{eq:wave cone}.

\subsection{Convex hull and geometric setup}
Given a cone $\Lambda$, we say that $K$ is convex with respect to $\Lambda$ if, for any two points 
$A, B \in K$ with $B-A\in \Lambda$, the whole segment $[A, B]$ belongs to $K$.
The $\Lambda$-convex hull of $K_{\rho,\chi}$ is the smallest $\Lambda$-convex set $K_{\rho,\chi}^{\Lambda}$ containing $K_{\rho,\chi}$, 
i.e. the set of states obtained by mixture of states of $K_{\rho,\chi}$ through oscillations in $\Lambda$-directions
(Gromov \cite{gr}, who works in the more general setting of jet bundles, calls this the $P$- convex hull).
The key point in Gromov's method of convex integration (which is a far reaching generalization of the work of Nash \cite{na}
and Kuiper \cite{ku} on isometric immersions) is that \eqref{eq:linear PDE} coupled with a pointwise nonlinear constraint of the form
$z\in K$ a.e. admits many interesting solutions provided that the $\Lambda$-convex hull of $K$, $K^{\Lambda}$, is sufficiently large.
In applications to elliptic and parabolic systems we always have $K^{\Lambda}=K$ so that Gromov's approach does not directly apply. For other
applications to partial differential equations it turns out that one can work with the $\Lambda$-convex hull defined by duality. More precisely, 
a point does not belong to the $\Lambda$- convex hull defined by duality if and only if there exists a $\Lambda$-convex function which separates it from $K$.
A crucial fact is that the second notion is much weaker. This surprising fact is illustrated in \cite{kms}.

In our case, the wave cone is quite large, therefore it is sufficient to consider the stronger notion of $\Lambda$-convex hull, indeed it coincides with 
the whole convex hull of $K_{\rho,\chi}$.

\begin{lemm} \label{lemma:convex hull}
For any $S\in\mathcal{S}^n$ let $\lambda_{max}(S)$ denote the
largest eigenvalue of $S$. For $(\rho, m, U)\in \R^+\times \R^n\times
\mathcal{S}^n_0$ let
\begin{equation} \label{eq:e}
e(\rho, m, U):= \lambda_{max}\left(\frac{m\otimes m}{\rho}-U\right).
\end{equation}
Then, for any given $\rho,\chi \in \R^+$, the following holds
\begin{itemize}
  \item[(i)] $e(\rho, \cdot, \cdot): \R^n\times \Sym_0^n  \rightarrow \R$ is convex;
  \item[(ii)] $ \frac{\abs{m}^2}{n\rho} \leq e(\rho, m, U)$, with equality if and only if $U=\frac{m\otimes m}{\rho}-\frac{\abs{m}^2}{n\rho} I_n$;
  \item[(iii)] ${\abs{U}}_\infty \leq (n-1) e(\rho, m, U)$, with ${\abs{U}}_\infty$ being the operator norm of the matrix;
  \item[(iv)] the $\frac{\chi}{n}$-sublevel set of $e$ defines the
convex hull of $K_{\rho,\chi}$, i.e.
\begin{align} \label{eq:convex hull}
K_{\rho,\chi}^{co}=\Bigg\{&(m, U, q)\in \R^n\times
\Sym_0^n\times \R^+: e(\rho, m, U)\leq \frac{\chi}{n},
\notag
\\  q&=p(\rho)+\frac{\chi}{n}\Bigg\}
\end{align}
and $K_{\rho,\chi}=K_{\rho,\chi}^{co} \cap \{\abs{m}^2=\rho \chi \}$.
\end{itemize}
\end{lemm}

For the proof of (i)-(iv) we refer the reader to the proof of Lemma
3.2 in \cite{dls2}: the arguments there can be easily adapted to
our case.

We observe that, for any $\rho,\chi\in \R^+$, the convex hull $K_{\rho,\chi}^{co}$ lives in the hyperplane H of $\R^n\times \Sym_0^n \times \R^+$ defined by 
$H:= \Big\{(m, U, q)\in \R^n\times \Sym_0^n\times \R^+: q=p(\rho)+\frac{\chi}{n}\Big\}.$
Therefore, the interior of $K_{\rho,\chi}^{co}$ as a subset of $\R^n\times \Sym_0^n \times \R^+$ is empty. This seems to prevent us 
from working in the classical framework of convex integration, but we can overcome this apparent obstacle.

For any $\rho,\chi\in \R^+$, we define the \textit{hyperinterior} of $K_{\rho,\chi}^{co}$, and we denote it with ``hint $K_{\rho,\chi}^{co}$'', as the following set
\begin{align} \label{eq:hyperinterior}
 \text{hint } K_{\rho,\chi}^{co}:=\Bigg\{&(m, U, q)\in \R^n\times
\Sym_0^n\times \R^+: e(\rho, m, U)< \frac{\chi}{n},
\notag
\\  q&=p(\rho)+\frac{\chi}{n}\Bigg\}.
\end{align}

In the framework of convex integration, the larger the $\Lambda$-convex hull of $K_{\rho,\chi}$ is,
the bigger the breathing space will be. How to ``quantify'' the meaning of a ``large'' $\Lambda$-convex hull in our context?
The previous definition provides an answer: the $\Lambda$-convex hull of $K_{\rho,\chi}$ will be ``large'' if its hyperinterior is nonempty.
The wave cone of the semi-stationary Euler isentropic system is wide enough to ensure that the $\Lambda$-convex hull of $K_{\rho,\chi}$ coincides with
the convex hull of $K_{\rho,\chi}$ and has a nonempty hyperinterior. As a consequence, we can construct irregular solutions oscillating
along any fixed direction. For our purposes, it will be convenient to restrict
to some \textit{special directions} in $\Lambda$, consisting of matrices of rank 2, which are not stationary
in time, but are associated with a constant pressure.

\begin{lemm}\label{lemma:special directions}
 Let $c,d\in\R^n$ with $\abs{c}=\abs{d}$ and $c\neq d$, and let $\rho\in\R^+$. \\Then
$\left(c-d, \frac{c\otimes c}{\rho}- \frac{d\otimes  d}{\rho}, 0 \right) \in
\Lambda$.
\end{lemm}

\begin{dimo}
Since the vector $\left(c+d, -\left(\frac{\abs{c}^2+c\cdot d}{\rho}\right)\right)$ is in the kernel of the matrix
$$
    C=\left( \begin{array}{cc}
    \frac{c\otimes c}{\rho}-\frac{d\otimes d}{\rho} & c-d \\
    c-d & 0\\
    \end{array} \right),
$$
$C$ has indeed determinant zero, hence $\left(c-d,\frac{c\otimes c}{\rho}- \frac{d\otimes  d}{\rho}, 0 \right)  \in \Lambda$.
\end{dimo}

Now, we introduce some important tools: they allow us to prove that $K_{\rho,\chi}^{\Lambda}=K_{\rho,\chi}^{co}$
is sufficiently large, thus providing us room to find many solutions for \eqref{eq:linearized system}-\eqref{eq:pointwise constraint}.

As first, we define the \textit{admissible segments} as segments in $\R^n\times \Sym_0^n\times \R^+$ whose
directions belong to the wave cone $\Lambda$ for the linear system of PDEs \eqref{eq:linearized system}
and are indeed \textit{special directions} in the sense specified by Lemma \ref{lemma:special directions}.

\begin{defin} \label{def:admissible segments}
 Given $\rho,\chi\in\R^+$ we call $\sigma$ an \textit{admissible segment} for $(\rho,\chi)$ if $\sigma$ is a line segment in
$\R^n\times \Sym_0^n\times \R^+$ satisfying the following conditions:
\begin{itemize}
 \item $\sigma$ is contained in the hyperinterior of $K_{\rho,\chi}^{co}$;
  \item $\sigma$ is parallel to $\left(c-d, \frac{c \otimes c}{\rho}-\frac{d\otimes d}{\rho}, 0\right)$ for some $c,d \in \R^n$ with
  $\abs{c}^2=\abs{d}^2=\rho\chi$ and $c\neq \pm d$.
\end{itemize}
\end{defin}

The admissible segments defined above correspond to suitable plane-wave solutions
of \eqref{eq:linearized system}. The following Lemma ensures that, for any $\rho,\chi\in \R^+$, the hyperinterior of $K_{\rho,\chi}^{co}$
is `` sufficiently round '' with respect to the special directions: given any point in the hyperinterior
of $K_{\rho,\chi}^{co}$, it can be seen as the midpoint of a sufficiently large admissible segment for $(\rho,\chi)$.

\begin{lemm}\label{lemma:Caratheodory lemma}
There exists a constant $F=F(n)>0$ such that for any $\rho,\chi\in \R^+$ and for any $z=(m, U, q)\in \text{hint } K_{\rho,\chi}^{co}$
there exists an admissible line segment for $(\rho,\chi)$
\begin{equation} \label{eq:admissible segment}
 \sigma=\left[ (m,U,q)- (\overline{m}, \overline{U}, 0), (m,U,q)+ (\overline{m}, \overline{U}, 0)\right]
\end{equation}
such that
$$\abs{\overline{m}}\geq \frac{ F }{\sqrt{\rho\chi}} \left(\rho\chi-\abs{m}^2 \right).$$
 \end{lemm}

The proof rests on a clever application of Carath\'{e}odory's theorem for convex sets and can be carried out, with minor modifications,
as in \cite{dls2} (cf. Lemma $6$ therein).

As an easy consequence of the previous Lemma, we can finally establish that the $\Lambda$-convex hull of $K_{\rho,\chi}$
coincides with $K_{\rho,\chi}^{co}$.

\begin{propos} \label{prop:lambda convex hull}
 For all given $\rho,\chi \in \R^+$, the $\Lambda$-convex hull of $K_{\rho,\chi}$ coincides with the convex hull of $K_{\rho,\chi}$.
\end{propos}

\begin{dimo}
Recall that, given $\rho,\chi\in\R^+$, we denote the $\Lambda$-convex hull of $K_{\rho,\chi}$ with $K_{\rho,\chi}^{\Lambda}$.
Of course $K_{\rho,\chi}^{\Lambda}\subset K_{\rho,\chi}^{co}$, hence we have to prove the opposite inclusion, i.e. $K_{\rho,\chi}^{co}\subset K_{\rho,\chi}^{\Lambda}$.
For every $z\in K_{\rho,\chi}^{co}$ we can follow the procedure in the proof of Lemma \ref{lemma:Caratheodory lemma} (cf. \cite{dls2}) and write it as
$z=\sum_j \lambda_j z_j$, with $(z_j)_{1\leq j\leq N+1}$ in $K_{\rho,\chi}$, $(\lambda_j)_{1\leq j\leq N+1}$ in $[0,1]$
and $\sum_j \lambda_j=1$. Again, we can assume that $\lambda_1=\max_j \lambda_j$. In case $\lambda_1=1$ then
$z=z_1\in K_{\rho,\chi} \subset K_{\rho,\chi}^{\Lambda}$ and we can already conclude.
Otherwise (i.e. when $\lambda_1\in (0,1)$) we can argue as in Lemma \ref{lemma:Caratheodory lemma} so to find an admissible
segment $\sigma$ for $(\rho,\chi)$ of the form \eqref{eq:admissible segment}.
Since we aim at writing $z$ as a $\Lambda$-barycenter of elements of $K_{\rho,\chi}$, we ``play'' with these admissible segments
by prolongations and iterative constructions until we get segments with extremes lying in $K_{\rho,\chi}$. More precisely:
we extend the segment $\sigma$ until we meet $\partial \text{ hint} K_{\rho,\chi}^{co}$ thus obtaining $z$ as the barycenter
of two points $(w_0, w_1)$ with $(w_0 - w_1)\in \Lambda$  and such that every $w_i=(m_i, U_i, q_i)$, $i=0,1$, satisfies
either $\abs{m_i}^2=\rho\chi$ or $\abs{m_i}^2<\rho\chi$ and $e(\rho, m_i, U_i)=\chi/ n$.

In the first case, $U_i -\left(\frac{m_i\otimes m_i}{\rho}-\frac{\abs{m_i}^2}{n\rho} I_n \right) \geq 0$, and since it is a null-trace-matrix
it is identically zero, whence $w_i\in K_{\rho,\chi}$ \Big(note that in the construction of $\sigma$ the $q$-direction remains unchanged, hence
$q_i= p(\rho)+\frac{\chi}{n}$\Big).

In the second case,  i.e. when $\abs{m_i}^2<\rho\chi$ and $e(\rho, m_i, U_i)=\chi/ n$, we apply again Lemma
\ref{lemma:Caratheodory lemma} and a limit procedure to express $w_i$ as barycentre of $(w_{i,0}, w_{i,1})$
with $(w_{i,0} - w_{i,1})\in \Lambda$  and such that every $w_{i,k}=(m_{i,k}, U_{i,k}, q_{i,k})$, $k=0,1$, will satisfy
either $\abs{m_{i,k}}^2=\rho\chi$ or $\lambda_2(\rho, m_{i,k}, U_{i,k})=e(\rho, m_{i,k}, U_{i,k})=\chi/ n$,
where $\lambda_1(\rho, m, U)\geq \lambda_2(\rho, m ,U) \geq......\geq \lambda_n(\rho, m, U)$ denote the ordered eigenvalues of the matrix
$\frac{m\otimes m}{\rho} -U$ (note that $\lambda_1(\rho, m , U)=e(\rho, m ,U)$).
Now, we iterate this procedure of constructing suitable admissible segments for $(\rho,\chi)$ until we have written $z$ as $\Lambda$-barycenter of
points $(m, U, q)$ satisfying either $\abs{m}^2=\rho\chi$ or $\lambda_n(\rho, m, U)=\chi/n$
and therefore all belonging to $K_{\rho,\chi}$ as desired.
\end{dimo}

\section{A criterion for the existence of infinitely many solutions}
\setcounter{equation}{0}
The following Proposition provides a criterion to recognize initial
data $m^0$ which allow for many weak admissible
solutions to \eqref{eq:isentropic Euler system}. Its proof relies deeply on the geometrical analysis carried out in Section 3. The underlying idea
comes from convex integration. The general principle of this method, developed
for partial differential equations by Gromov \cite{gr} and for ordinary differential equations
by Filippov \cite{fi}, consists in the following steps:
given a nonlinear equation $\mathcal{E}(z)$,
\begin{itemize}
 \item (i) we rewrite it as $(\mathcal{L}(z) \wedge z\in K)$ where $\mathcal{L}$ is a linear equation;
\item (ii) we introduce a strict subsolution $z_0$ of the system, i.e. satisfying a relaxed system $(\mathcal{L}(z_0) \wedge z\in \mathcal{U})$;
\item (iii) we construct a sequence $(z_k)_{k\in \mathbb{N}}$ approaching $K$ but staying in $\mathcal{U}$;
\item (iv) we pass to the limit, possibly modifying the sequence $(z_k)$ in order to ensure a suitable convergence.
\end{itemize}

Step (i) has already been done in Section 3.1. The choice of $z_0$ will be specified in Sections 7-8. Here, we define the notion of subsolution for an appropriate 
set $\mathcal{U}$, we construct an improving sequence and we pass to the limit. The way how we construct the approximating sequence will be
described in Section 6 using some tools from Section 5.

One crucial step in convex integration is the passage from open sets $K$ to general sets. This can be done in different ways, e.g.
by the Baire category theorem (cf. \cite{ox}), a refinement of it using Baire-1 functions or the Banach-Mazur game \cite{ki} or by direct construction \cite{sy}.
Whatever approach one uses the basic theme is the same: at each step of the construction one adds a highly oscillatory correction whose frequency is much larger and whose
amplitude is much smaller than those of the previous corrections.

In this section, we achieve our goals following some Baire category arguments as in \cite{dls1}: they
are morally close to the methods developed by Bressan and Flores in
\cite{bf} and by Kirchheim in \cite{ki}. 

In our framework the initial data will be
constructed starting from solutions to the convexified (or relaxed) problem
associated to \eqref{eq:semistationary system}, i.e. solutions to
the linearized system \eqref{eq:linearized system} satisfying a
``relaxed'' nonlinear constraint \eqref{eq:pointwise
constraint} (i.e. belonging to the hyperinterior of the convex hull of the ``constraint set''),
which we will call \textit{subsolutions}.

As in \cite{dls1}, our application shows that the Baire theory is comparable in terms of results to the method of
convex integration and they have many similarities: they are both
based on an approximation approach to tackle problems while the
difference lies only in the limit arguments, i.e. on the way the exact solution is obtained from better and better approximate ones. 
These similarities are
clarified by Kirchheim in \cite{ki}, where the continuity
points of a first category Baire function are considered; a
comparison between the two methods is drawn by Sychev in
\cite{sy}. 

Here, the topological reasoning of Baire
theory is preferred to the iteration technique of convex
integration, since the first has the advantage to provide us
directly with infinitely many different solutions.

\begin{propos}\label{prop:criterion}
Let $\rho_0\in C^1_p(Q;\R^+)$ be a given density function and let $T$ be any finite positive time. Assume there exist $(m_0,
U_0, q_0)$ continuous space-periodic solutions of \eqref{eq:linearized system} on
$\R^n\times ]0,T[$ with 
\begin{equation} \label{eq:hypothesis 1 criterion}
m_0 \in C([0,T];H_w(Q)), 
\end{equation}
and a function $\chi\in C^\infty([0,T]; \R^+)$ such that
\begin{align}\label{eq:hypothesis 3 criterion}
&e(\rho_0(x),m_0(x,t),U_0(x,t))<\frac{\chi(t)}{n} \qquad\text{for all } (x,t)\in \R^n\times ]0,T[,\\
&q_0(x,t)=p(\rho_0(x))+\frac{\chi(t)}{n} \qquad\text{for all } (x,t)\in \R^n\times ]0,T[.
\end{align}
Then there exist infinitely many weak
solutions $(\rho,m)$ of the system \eqref{eq:semistationary system} in $\R^n\times [0,T[$
with density $\rho(x)=\rho_0(x)$ and such that
\begin{align}
&m\in C([0,T];H_w(Q)), \label{eq:thesis 1 criterion}\\
&m(\cdot,t)=m_0(\cdot, t) \qquad \text{for }t=0,T \text{ and for a.e. } x\in \R^n,
\label{eq:thesis 2 criterion} \\
&\abs{m(x,t)}^2=\rho_0(x)\chi(t)\qquad \text{for
a.e. } (x,t)\in \R^n\times ]0,T[ . \label{eq:thesis
3 criterion}
\end{align}
\end{propos}

\subsection{The space of subsolutions}
We define the space of \textit{subsolutions} as follows. Let $\rho_0$ and $\chi$ be given as in the assumptions of
Proposition \ref{prop:criterion}. Let $m_0$
be a vector field as in Proposition \ref{prop:criterion} with
associated modified pressure $q_0$ and consider space-periodic momentum fields $m:\R^n\times
[0,T]\rightarrow \R^n$ which satisfy
\begin{equation} \label{eq:subsolution1}
\div m=0,
\end{equation}
the initial and boundary conditions
\begin{align}
&m(x,0)=m_0(x,0),\label{eq:subsolution2}\\
&m(x,T)=m_0(x,T),\label{eq:subsolution3}\\
\end{align}
and such that there exists a continuous space-periodic matrix field $U:\R^n\times
]0,T[\rightarrow \Sym_0^n$ with 
\begin{equation}\label{eq:subsolution5}
\begin{aligned} 
&e(\rho_0(x),m(x,t), U(x,t))<\frac{\chi}{n}\quad \text{for all } (x,t)\in\R^n\times ]0,T[,\\
&\dt m+\div_x U+\nabla_x q_0 = 0\quad \text{in }\R^n\times [0,T].
\end{aligned}
\end{equation}

\begin{defin}\label{def:subsolutions}
Let $X_0$ be the set of such linear momentum fields, i.e.
\begin{align} \label{eq:X_0}
X_0=\Bigg\{&m \in C^0(]0,T[; C^0_p(Q))\cap C([0,T];H_w(Q)):\notag \\
&\eqref{eq:subsolution1}-\eqref{eq:subsolution5}\text{ are
satisfied}\Bigg\}
\end{align}
and let $X$ be the closure of $X_0$ in $C([0,T]; H_w(Q)$. Then $X_0$
will be the space of \textit{strict subsolutions}.
\end{defin}

As $\rho_0$ is continuous and periodic on $\R^n$ and $\chi$ is smooth on $[0,T]$, there exists a constant $G$ such
that $\chi(t) \int_Q \rho_0(x) dx \leq G$ for all $t\in [0,T]$. Since
for any $m\in X_0$ with associated matrix field $U$ we have that
(see Lemma \ref{lemma:convex hull}- (ii))
\begin{align*}
\int_Q \abs{m(x,t)}^2 dx &\leq \int_Q n \rho_0(x) e(\rho_0(x),
m(x,t), U(x,t)) dx \\
&< \chi(t) \int_Q \rho_0(x) dx \text{ for all }
t\in [0,T],
\end{align*}
we can observe that $X_0$ consists of functions
$m:[0,T]\rightarrow H(Q)$ taking values in a bounded subset
$B$ of $H(Q)$. Without loss of generality, we can assume that
$B$ is weakly closed. Then, $B$ in its weak topology is metrizable
and, if we let $d_B$ be a metric on $B$ inducing the weak topology,
we have that $(B, d_B)$ is a compact metric space. Moreover, we can
define on $Y:=C([0,T], (B, d_B))$ a metric $d$ naturally induced by
$d_B$ via
\begin{equation} \label{eq:metric d}
d(f_1, f_2):= \max _{t\in [0,T]} d_B (f_1(\cdot, t), f_2 (\cdot,
t)).
\end{equation}
Note that the topology induced on $Y$ by $d$ is equivalent to the
topology of $Y$ as a subset of $C([0,T];H_w)$. In addition,
the space $(Y,d)$ is complete. Finally, $X$
is the closure in $(Y,d)$ of $X_0$ and hence $(X,d)$ is as well a
complete metric space.

\begin{lemm}\label{lemma:set of subsolutions}
If $m\in X$ is such that $\abs{m(x,t)}^2=\rho_0(x) \chi(t)$ for almost every $(x,t)\in \R^n \times ]0,T[$, then the pair $(\rho_0,m)$
is a weak solution of \eqref{eq:semistationary system} in $\R^n\times [0,T[ $ satisfying \eqref{eq:thesis 1 criterion}-\eqref{eq:thesis 2 criterion}-\eqref{eq:thesis 3 criterion}.
\end{lemm}

\begin{dimo}
Let $m\in X$ be such that $\abs{m(x,t)}^2=\rho_0(x) \chi(t)$ for almost every $(x,t)\in \R^n \times ]0,T[$.
By density of $X_0$, there exists a sequence $\{m_k\}\subset X_0$ such that
$m_k\stackrel{d}{\rightarrow}m$ in $X$. For any $m_k \in X_0$ let $U_k$ be the associated smooth matrix field
enjoying \eqref{eq:subsolution5}. Thanks to Lemma \ref{lemma:convex hull} (iii) and \eqref{eq:subsolution5},
the following pointwise estimate holds for the sequence $\{U_k\}$
$$\abs{U_k}_\infty \leq (n-1) e(\rho_0,m_k,U_k)< \frac{(n-1)-\chi}{n}.$$
As a consequence, $\{U_k\}$ is uniformly bounded in $L^\infty([0,T]; L^\infty_p(Q))$; by possibly extracting a subsequence, we
have that
$$U_k\stackrel{\ast}{\rightharpoonup} U \text{ in } L^\infty([0,T]; L^\infty_p(Q)).$$
Note that $\overline{\text{hint} K_{\rho_0,\chi}^{co}}=K_{\rho_0,\chi}^{co} $ is a convex and compact set by Lemma \ref{lemma:convex hull}-(i)-(ii)-(iii). 
Hence, $m\in X$ with associated matrix field $U$ solves \eqref{eq:linearized system} on $\R^n\times [0,T]$ for $q=q_0$ and $(m,U, q_0)$ takes values in $K_{\rho_0,\chi}^{co}$ almost everywhere. If, in addition, 
$\abs{m(x,t)}^2=\rho_0(x) \chi(t)$, then $(m,U, q_0)(x,t)\in K_{\rho,\chi}$ a.e. in $\R^n\times [0,T]$ (cf. Lemma \ref{lemma:convex hull}-(iv)). Lemma \ref{lemma:differential inclusion} 
allows us to conclude that  $(\rho_0,m)$
is a weak solution of \eqref{eq:semistationary system} in $\R^n\times [0,T[$.
Finally, since $m_k \rightarrow m$ in $C([0,T]; H_w(Q))$ and $\abs{m(x,t)}^2=\rho_0(x) \chi(t)$ for almost every $(x,t)\in \R^n \times ]0,T[$, 
we see that $m$ satisfies also \eqref{eq:thesis 1 criterion}-\eqref{eq:thesis 2 criterion}-\eqref{eq:thesis 3 criterion}.
\end{dimo}

Now, we will argue as in \cite{dls1} exploiting Baire category techniques to combine weak and strong convergence (see also \cite{ki}).

\begin{lemm}\label{lemma:Baire-1 map}
 The identity map $I:(X, d)\rightarrow L^2 ([0,T]; H(Q))$ defined by $m\rightarrow m$ is a Baire-1 map, and therefore the set of points of continuity is residual in
$(X,d)$.
\end{lemm}

\begin{dimo}
 Let $\phi_r(x,t)=r^{-(n+1)} \phi(rx, rt)$ be any regular spacetime convolution kernel.
For each fixed $m\in X$, we have
$$\phi_r \ast m \rightarrow m \quad \text{strongly in } L^2(H) \text{ as } r\rightarrow 0.$$
On the other hand, for each $r>0$ and $m_k\in X$,
$$m_k \stackrel{d}{\rightarrow} m \text{ implies } \phi_r\ast m_k \rightarrow \phi_r\ast m \text{ in } L^2(H).$$
Therefore, each map $I_r:(X,d)\rightarrow L^2 (H)$, $m\rightarrow \phi_r\ast m$ is continuous, and $I(m)= \lim_{r\rightarrow 0} I_r (m)$ for all $m\in X$.
This shows that $I:(X,d)\rightarrow L^2 (H)$ is a pointwise limit of continuous maps; hence it is a Baire-1 map. As a consequence, the set of points of continuity of $I$
is residual in $(X,d)$ (cf. \cite{ox}).
\end{dimo}

\subsection{Proof of Proposition \ref{prop:criterion}}
We aim to show that all points of continuity of the identity map correspond to solutions of \eqref{eq:semistationary system} enjoying the requirements 
of Proposition \ref{prop:criterion}: Lemma \ref{lemma:Baire-1 map} will then allow us to prove Proposition \ref{prop:criterion} once we know that the cardinality of 
$X$ is infinite.
In light of Lemma \ref{lemma:set of subsolutions}, for our purposes it suffices to prove the following claim:
\begin{claim} 
If $m\in X$ is a point of continuity of $I$, then 
\begin{equation} \label{eq:claim}
 \abs{m(x,t)}^2=\rho_0(x) \chi(t) \text{ for almost every }(x,t)\in \R^n \times ]0,T[.
\end{equation}
\end{claim}

Note that proving \eqref{eq:claim} is equivalent to prove that $\norm{m}_{L^2(Q\times [0,T])}= \left(\int_{Q}\int_0^T \rho_0(x) \chi(t)  dt dx\right)^{1/2}$, since for any $m\in X$ we have 
$\abs{m(x,t)}^2\leq \rho_0(x) \chi(t)$ for almost all $(x,t)\in \R^n\times [0,T]$.
Thanks to this remark, the claim is reduced to the following lemma (cf. Lemma 4.6 in \cite{dls1}), which provides a strategy to move towards the boundary of $X_0$: given $m\in X_0$,
we will be able to approach it with a sequence inside $X_0$ but closer than $m$ to the boundary of $X_0$. 

\begin{lemm} \label{lemma:improvement lemma}
Let $\rho_0, \chi$ be given functions as in Proposition \ref{prop:criterion}.
Then, there exists a constant $\beta=\beta(n)$ such that, given $m\in X_0$, there exists a sequence $\{m_k\}\subset X_0$ with the following properties
\begin{align} \label{eq:improvement}
 \norm{m_k}_{L^2(Q\times [0,T])}^2\geq  &\norm{m}_{L^2(Q\times [0,T])}^2 \notag \\
&+\beta \left( \int_{Q}\int_0^T \rho_0(x) \chi(t)  dt dx - \norm{m}_{L^2(Q\times [0,T])}^2\right)^2
\end{align}
and
\begin{equation} \label{eq:convergence of improvement}
 m_k \rightarrow m \text{ in } C([0,T], H_w(Q)).
\end{equation}
\end{lemm}

The proof is postponed to Section \ref{sec:improvement step}. Let us show how Lemma \ref{lemma:improvement lemma} implies the claim.
As in the claim, assume that $m\in X$ is a point of continuity of the identity map $I$. Let $\{m_k\}\subset X_0$ be a fixed sequence that converges to $m$
in $C([0,T], H_w(Q))$. Using Lemma \ref{lemma:improvement lemma} and a standard diagonal argument, we can find a second sequence $\{\widetilde{m}_k\}$
yet converging to $m$ in $X$ and satisfying
\begin{align*}
 \liminf _{k\rightarrow \infty} &\norm{\widetilde{m}_k}_{L^2(Q\times [0,T])}^2\geq  \liminf_{k\rightarrow\infty} \Bigg( \norm{m_k}_{L^2(Q\times [0,T])}^2 \notag \\
&+\beta \left( \int_{Q}\int_0^T \rho_0(x) \chi(t)  dt dx - \norm{m_k}_{L^2(Q\times [0,T])}^2\right)^2 \Bigg).
\end{align*}
According to the hypothesis, $I$ is continuous at $m$, therefore both $m_k$ and $\widetilde{m}_k$ converge strongly to $m$ and
\begin{align*}
 \norm{m}_{L^2(Q\times [0,T])}^2 &\geq \norm{m}_{L^2(Q\times [0,T])}^2 \notag \\  &+\beta \left( \int_{Q}\int_0^T \rho_0(x) \chi(t)  dt dx - \norm{m}_{L^2(Q\times [0,T])}^2\right)^2.
\end{align*}
Hence $\norm{m}_{L^2(Q\times [0,T])}= \left(\int_{Q}\int_0^T \rho_0(x) \chi(t)  dt dx\right)^{1/2}$ and the claim holds true.
Finally, since the assumptions of Proposition  \ref{prop:criterion} ensure that $X_0$ is nonempty, by Lemma \ref{lemma:improvement lemma} we can see that the cardinality of $X$
is infinite whence the cardinality of any residual set in $X$ is infinite. In particular, the set of continuity points of $I$ is infinite: this and the claim conclude the proof of
Proposition \ref{prop:criterion}.

\section{Localized oscillating solutions}
\setcounter{equation}{0}
The wild solutions are made by adding one dimensional oscillating functions 
in different directions $\lambda\in\Lambda$. For that it is needed to localize the waves.
More precisely, the proof of Lemma \ref{lemma:improvement lemma} relies
on the construction of solutions to the linear system
\eqref{eq:linearized system}, localized in space-time and
oscillating between two states in $K_{\rho_0,\chi}^{co}$ along a given special direction
$\lambda\in \Lambda$. Aiming at compactly supported solutions, one
faces the problem of localizing vector valued functions: this is
bypassed thanks to the construction of a ``localizing'' potential
for the conservation laws \eqref{eq:linearized system}. This
approach is inherited from \cite{dls2}. As in
\cite{dls1} it could be realized for every $\lambda\in
\Lambda$, but in our framework it is convenient to restrict only to
\textit{special} $\Lambda$-directions (cf. \cite{dls2}): this restriction will allow us to localize the oscillations at
constant pressure.

Why oscillations at constant pressure are meaningful for us and
needed in the proof of Lemma \ref{lemma:improvement lemma}?

Owing to Section 3, in the variables $y=(x,t)\in \R^{n+1}$, the
system \eqref{eq:linearized system} is equivalent to $\div_y M =0$,
where $M\in \Sym^{n+1}$ is defined via the linear map
\begin{equation} \label{eq:linear map}
    \R^n \times \Sym_0^n \times \R \ni (m, U, q) \longmapsto M=\left( \begin{array}{cc}
    U+q I_n & m \\
    m & 0\\
    \end{array} \right).
\end{equation}
More precisely, this map builds an identification between the set of solutions $(m,U,q)$ to \eqref{eq:linearized system}
and the set of symmetric $(n+1)\times(n+1)$ matrices $M$ with $M_{(n+1)(n+1)}=0$ and $\tr (M)=q$.

Therefore, solutions of \eqref{eq:linearized system} with $q\equiv0$
correspond to matrix fields
$M:\R^{n+1}\rightarrow\R^{(n+1)\times(n+1)}$ such that
\begin{equation} \label{eq:linearized system in matrix form}
\div_y M = 0, \quad M^{\text{T}}=M, \quad M_{(n+1)(n+1)}=0, \quad \tr (M)=0.
\end{equation}
Moreover, given a density $\rho$ and two states $(c, U_c, q_c),(d,
U_d, q_d)\in K_{\rho}$ with non collinear momentum vector fields $c$
and $d$ having same magnitude ($\abs{c}=\abs{d}$), and hence same
pressure ($q_c=q_d$), then the corresponding matrices $M_c$ and
$M_d$ have the following form
\begin{equation*}
 M_c=\left( \begin{array}{cc}
    \frac{c\otimes c}{\rho}+p (\rho)I_n & c \\
    c & 0\\
    \end{array} \right) \text{ and }
 M_d=\left( \begin{array}{cc}
    \frac{d\otimes d}{\rho}+p (\rho)I_n & d \\
    d & 0\\
    \end{array} \right)
\end{equation*}
and satisfy
\begin{equation*}
 M_c-M_d=\left( \begin{array}{cc}
    \frac{c\otimes c}{\rho}-\frac{d\otimes d}{\rho} & c-d \\
    c-d & 0\\
    \end{array} \right).
\end{equation*}
Finally note that $\tr(M_c-M_d)=0$ and $M_c-M_d\in \Lambda$
corresponds to a \textit{special direction}.

The following Proposition provides a potential for solutions of
\eqref{eq:linearized system} oscillating between two states $M_c$
and $M_d$ at constant pressure. It is an easy adaptation to our
framework of Proposition 4 in \cite{dls2}. 

\begin{propos} \label{prop:potential}
 Let $c,d\in\R^n$ such that $\abs{c}=\abs{d}$ and $c\neq d$. Let also $\rho\in \R$. Then
 there exists a matrix-valued, constant coefficient, homogeneous
 linear differential operator of order $3$
 $$ A(\partial): C_c^\infty (\R^{n+1})\rightarrow C_c^\infty
 (\R^{n+1}; \R^{(n+1)\times(n+1)})$$
 such that $M=A(\partial)\phi$ satisfies \eqref{eq:linearized system
 in matrix form} for all $\phi \in C_c^\infty (\R^{n+1})$. Moreover
 there exists $\eta\in \R^{n+1}$ such that
 \begin{itemize}
 \item  $\eta$ is not parallel to $e_{n+1}$;
 \item  if $\phi(y)=\psi(y\cdot \eta)$, then
 $$A(\partial )\phi (y)= (M_c-M_d) \phi '''(y \cdot \eta).$$
\end{itemize}
\end{propos}

We also report Lemma 7 from \cite{dls2}: it ensures that the oscillations of the planewaves generated in proposition
\ref{prop:potential} have a certain size in terms of an appropriate norm-type-functional.

\begin{lemm} \label{lemma:potential}
 Let $\eta\in \R^{n+1}$ be a vector which is not parallel to $e_{n+1}$. Then for any bounded open set $B\subset \R^n$
$$\lim_{N\rightarrow \infty} \int_B \sin^2(N \eta\cdot (x,t)) dx=\frac{1}{2} \abs{B}$$
uniformly in $t\in\R$.
\end{lemm}

For the proof we refer the reader to \cite{dls2}.

\section{The improvement step} \label{sec:improvement step}
\setcounter{equation}{0}
We are now about to prove one of the cornerstones of the costruction.
Before moving forward, let us resume the plan. 
We have already identified a relaxed problem by introducing \textit{subsolutions}. Then, we 
have proved a sort of ``\textit{h}-principle'' (even if there is no homotopy here)
according to which, the space of \textit{subsolutions} can be ``reduced'' to the space of solutions or, equivalently, the typical (in Baire's sense) \textit{subsolution} is a solution.
Once assumed that a \textit{subsolution} exists, the proof of our ``\textit{h}-principle'' builds upon Lemma \ref{lemma:improvement lemma} 
combined with Baire category arguments.
Indeed, we could also prove Proposition \ref{prop:criterion} by applying iteratively Lemma \ref{lemma:improvement lemma} and thus constructing
a converging sequence of subsolutions approaching $K_{\rho,\chi}$: this would 
correspond to the convex integration approach.
So two steps are left in order to conclude our argument: showing the existence of a ``starting'' \textit{subsolution} and prove 
Lemma \ref{lemma:improvement lemma}.

This section is devoted to the second task, the proof of Lemma \ref{lemma:improvement lemma}, 
while in next section we will exhibit a ``concrete''
\textit{subsolution}.

What follows will be quite technical, therefore we first would like to recall the plan: we will add fast oscillations in 
allowed directions so to let $\abs{m}^2$ increase in average.
The proof is inspired by \cite{dls1}-\cite{dls2}.

\begin{dimo} [Proof of Lemma \ref{lemma:improvement lemma}]
Let us fix the domain $\Omega:=Q\times[0,T]$.
We look for a sequence $\{m_k\}\subset X_0$, with associated matrix fields $\{U_k\}$, which improves $m$ in the sense
of \eqref{eq:improvement} and has the form
\begin{equation} \label{eq:improving sequence}
(m_k, U_k)=(m, U)+\sum_{j} (\widetilde{m}_{k,j}, \widetilde{U}_{k,j})
\end{equation}
where every $z_{k,j}=(\widetilde{m}_{k,j}, \widetilde{U}_{k,j})$ is compactly supported in some suitable ball $B_{k,j}(x_{k,j}, t_{k,j})\subset\Omega$.
We proceed as follows.

\textit{Step 1.} 
Let $m\in X_0$ with associated matrix field $U$. By Lemma \ref{lemma:Caratheodory lemma}, for any $(x,t)\in \Omega$ 
we can find a line segment $\sigma_{(x,t)}:=[(m(x,t), U(x,t), q_0(x))- \lambda_{(x,t)}, (m(x,t), U(x,t), q_0(x))+ \lambda_{(x,t)}]$
admissible for $(\rho_0(x), \chi(t))$ and with direction $$\lambda_{(x,t)} =(\overline{m}(x,t), \overline{U}(x,t), 0)$$ such that
\begin{equation} \label{eq:Caratheodory estimate}
\abs{\overline{m}(x,t)} \geq \frac{F}{\sqrt{\rho_0(x) \chi(t)}} \left(\rho_0(x) \chi(t)-\abs{m(x,t)}^2\right).
\end{equation}
Since $z:=(m,U)$ and $K_{\rho_0,\chi}^{co}$ are uniformly continuous in $(x,t)$, there exists an $\varepsilon>0$ such that for any $(x,t), (x_0,t_0)\in \Omega$
with $\abs{x-x_0}+\abs{t-t_0}<\varepsilon$, we have
\begin{equation} \label{eq:perturbation}
(z(x,t), q_0(x)) \pm (\overline{m}(x_0, t_0), \overline{U}(x_0, t_0), 0) \subset \text{hint} K_{\rho_0,\chi}^{co}.
\end{equation}

\textit{Step 2.} 
Fix $(x_0,t_0)\in \Omega$ for the moment. Now, let $0\leq\phi_{r_0}\leq 1$ be a smooth cutoff function on $\Omega$ with support contained 
in a ball $B_{r_0}(x_0, t_0)\subset\Omega$ for some $r_0>0$, identically 1 on $B_{r_0/2}(x_0, t_0)$ and strictly less than 1 outside.
Thanks to Proposition \ref{prop:potential} and the identification $(m, U,q)\rightarrow M$, for the admissible line segment $\sigma_{(x_0,t_0)}$,
there exist an operator $A_0$ and a direction $\eta_0 \in \R^{n+1}$ not parallel to $e_{n+1}$, such that for any $k\in\N$
\begin{equation*}
 A_0 \left(\frac{\cos (k \eta_0 \cdot (x,t))}{k^3}\right)=\lambda_{(x_0,t_0)} \sin (k \eta_0 \cdot (x,t)),
\end{equation*}
and such that the pair $(\widetilde{m}_{k,0}, \widetilde{U}_{k,0})$ defined by 
\begin{equation*}
 (\widetilde{m}_{k,0}, \widetilde{U}_{k,0})(x,t):= A_0 \left[\phi_{r_0} (x,t) k^{-3} \cos(k \eta_0 \cdot (x,t))\right]
\end{equation*}
satisfies \eqref{eq:linearized system} with $q \equiv 0 $.
Note that $(\widetilde{m}_{k,0}, \widetilde{U}_{k,0})$ is supported in the ball $B_{r_0}(x_0, t_0)$ and that
\begin{align} \label{eq:estimate in L infinite}
& \norm{(\widetilde{m}_{k,0}, \widetilde{U}_{k,0})-\phi_{r_0} \left(\overline{m}(x_0,t_0), \overline{U}(x_0,t_0)\right) \sin (k\eta_0 \cdot (x,t))}_{\infty} \notag \\
& \leq const \left( A_0, \eta_0, \norm{\phi_0}_{C^3}\right) \frac{1}{k}
\end{align}
since $A_0$ is a linear differential operator of homogeneous degree $3$. 
Furthermore, for all $(x,t)\in B_{r_0/2}(x_0, t_0)$, we have
\begin{equation*}
 \abs{\widetilde{m}_{k,0}(x,t)}^2=\abs{\overline{m}(x_0,t_0)}^2 \sin^2(k \eta_0  \cdot (x,t)).
\end{equation*}
Since $\eta_0 \in \R^{n+1}$ is not parallel to $e_{n+1}$, from Lemma \ref{lemma:potential} we can see that
\begin{equation*}
 \lim_{k\rightarrow \infty} \int_{B_{r_0/2}(x_0, t_0)} \abs{\widetilde{m}_{k,0}(x,t)}^2 dx=\frac{1}{2} \int_{B_{r_0/2}(x_0, t_0)} \abs{\overline{m}(x_0,t_0)}^2 dx
\end{equation*}
uniformly in $t$. In particular, using \eqref{eq:Caratheodory estimate}, we obtain
\begin{align}  \label{eq:eq1}
 \lim_{k\rightarrow \infty} &\int_{B_{r_0/2}(x_0, t_0)} \abs{\widetilde{m}_{k,0}(x,t)}^2 dx dt \geq \notag  \\
&\frac{F^2}{2\rho_0(x_0)\chi(t_0)} \left(\rho_0(x_0)\chi(t_0)) -\abs{m(x_0, t_0)}^2\right)^2 \abs{B_{r_0/2}(x_0, t_0)}.
\end{align}

\textit{Step 3.}
Next, observe that since $m$ is uniformly continuous, there exists an $\bar{r}>0$ such that for any $r< \bar{r}$ there exists a finite family of pairwise disjoint balls 
$B_{r_j}(x_j,t_j)\subset \Omega$ with $r_j<\bar{r}$ such that
\begin{align} \label{eq:eq2}
 \int_\Omega &\left( \rho_0 (x) \chi(t)- \abs{m(x,t)}^2 \right)^2 dx dt \leq \notag \\
&2 \sum_j \left( \rho_0(x_j) \chi(t_j) - \abs{m(x_j,t_j)}^2 \right)^2 \abs{B_{r_j}(x_j,t_j)}.
\end{align}
Fix $s>0$ with $s < \min\{\bar{r}, \varepsilon\}$ and choose a finite family of pairwise disjoint balls $B_{r_{j}}(x_{j},t_{j})\subset \Omega$ 
with radii $r_{j} < s$ such that \eqref{eq:eq2} holds.
In each ball $B_{2r_{j}}(x_{j},t_{j})$ we apply the construction of \textit{Step 2} to obtain, for every $k\in \N$, a pair $(\widetilde{m}_{k,j}, \widetilde{U}_{k,j})$.

\textit{Final step.}
Letting $(m_k, U_k)$ be as in \eqref{eq:improving sequence}, we observe that the sum therein consists of finitely many terms.
Therefore from \eqref{eq:perturbation} and \eqref{eq:estimate in L infinite} we deduce that there exists $k_0\in \N$ such that
\begin{equation} \label{eq:condition 1}
m_k\in X_0 \text{ for all } k\geq k_0. 
\end{equation}
Moreover, owing to \eqref{eq:eq1} and \eqref{eq:eq2} we can write
\begin{align} \label{eq:eq3}
\lim_{k\rightarrow \infty} &\int_\Omega \abs{m_k(x,t) -m (x,t)}^2 dx dt= \sum_j \lim_{k\rightarrow \infty} \int_\Omega \abs{\widetilde{m}_{k,j}(x,t)}^2 dx dt \notag \\
&\geq \sum_j \frac{F^2}{2\rho_0(x_j)\chi(t_j)} \left(\rho_0(x_j)\chi(t_j)) -\abs{m(x_j, t_j)}^2\right)^2 \abs{B_{r_j}(x_j, t_j)}\notag  \\
&\geq C \int_\Omega \left( \rho_0 (x) \chi(t)- \abs{m(x,t)}^2 \right)^2 dx dt.
\end{align}
Since $m_k\stackrel{d}{\rightarrow} m$, due to \eqref{eq:eq3} we have
\begin{align} \label{eq:eq4}
\liminf_{k\rightarrow \infty} \norm{m_k}^2_{L^2(\Omega)} &= \norm{m}_2^2 +\liminf_{k\rightarrow \infty} \norm{m_k-m}_2^2 \notag\\
&\geq \norm{m}_2^2 + C \int_\Omega \left( \rho_0 (x) \chi(t)- \abs{m(x,t)}^2 \right)^2 dx dt,
\end{align}
which gives \eqref{eq:improvement} with $\beta=\beta(n)=\beta(F(n))$.
\end{dimo}

\section{Construction of suitable initial data}
\setcounter{equation}{0}
In this section we show the existence of a \textit{subsolution} in the sense of Definition
\ref{def:subsolutions}. 
Since the subsolution we aim to construct has to be space-periodic, it will be enough to work on the building brick
$Q$ and then extend the costruction periodically to $\R^n$. 

The idea to work in the space-periodic setting has been recently adopted by Wiedemann \cite{wi}
in order to construct global solutions to the incompressible Euler equations.

\begin{propos} \label{prop:starting subsolution}
Let $\rho_0\in C^1_p(Q;\R^+)$ be a given density function as in Proposition \ref{prop:criterion} and let $T$ be any given positive time.
Then, there exist a smooth function $\widetilde{\chi}:\R\rightarrow \R^+$, 
a continuous periodic matrix field $\widetilde{U}: \R^n \rightarrow \Sym^n_0$ and a function $\widetilde{q}\in C^1(\R; C^1_p(\R^n))$ such that
\begin{equation} \label{eq:linearized equations1}
    \div_x \widetilde{U}+\nabla_x \widetilde{q}=0 \qquad \text{on} \quad \R^n\times \R
\end{equation}
and
\begin{align}
&e(\rho_0(x),0,\widetilde{U}(x))<\frac{\widetilde{\chi}(t)}{n}\text{
for all } (x,t)\in \R^n\times [0,T[ \label{eq:subsolutions
constraint1}\\
&\widetilde{q}(x,t)=p(\rho_0(x))+\frac{\widetilde{\chi}(t)}{n} \text{ for all } x\in \R^n\times \R. \label{eq:q
condition}
\end{align}
\end{propos}

\begin{dimo}[Proposition \ref{prop:starting subsolution}]
Let us define $\widetilde{U}$ componentwise by its Fourier transform as follows:
\begin{align} \label{eq:fourier transform of u}
 &\widehat{\widetilde{U}}_{ij}(k):= \Bigg( \frac{n k_i k_j}{(n-1)\abs{k}^2} \Bigg) \widehat{p(\rho_0(k))} \text{ if } i\neq j, \notag \\
&\widehat{\widetilde{U}}_{ii}(k):=\Bigg( \frac{n k_i^2 -\abs{k}^2}{(n-1) \abs{k}^2} \Bigg) \widehat{p(\rho_0(k))}.
\end{align}
for every $k\neq 0$, and $\widehat{\widetilde{U}}(0)=0$.
Clearly $\widehat{\widetilde{U}}_{ij}$ thus defined is symmetric and trace-free.
Moreover, since $p(\rho_0)\in C^1_p(\R^n)$, standard elliptic regularity arguments allow us to conclude 
that $\widetilde{U}$ is a continuous periodic matrix field.
Next, notice that 
\begin{equation} \label{eq:max eigenvalue}
 \norm{e(\rho_0(x),0,\widetilde{U}(x))}_\infty= \norm{\lambda_{max}(-\widetilde{U})}_\infty= \widetilde{\lambda}
\end{equation}
for some positive constant $\widetilde{\lambda}$.
Therefore, we can choose any smooth function $\widetilde{\chi}$ on $\R$ such that $\widetilde{\chi} > n\widetilde{\lambda}$ on $[0,T]$ in order to ensure 
\eqref{eq:subsolutions constraint1}. 
Now, let $\widetilde{q}$ be defined exactly as in \eqref{eq:q condition} for the choice of $\widetilde{\chi}$ just done.
It remains to show that \eqref{eq:linearized equations1} holds.
In light of \eqref{eq:q condition}, we can write equation \eqref{eq:linearized equations1} in Fourier space as
\begin{equation} \label{eq:eq in fourier space}
 \sum_{j=1}^n k_j \widehat{\widetilde{U}}_{ij}= k_i \widehat{p(\rho_0)}
\end{equation}
for $k\in\Z^n$.
It is easy to check that $\widehat{\widetilde{U}}$ as defined by \eqref{eq:fourier transform of u} solves
\eqref{eq:eq in fourier space} and hence $\widetilde{U}$ and $\widetilde{q}$ satisfy \eqref{eq:linearized equations1}

\end{dimo}

\begin{obs} 
We note that the H\"{o}lder continuity of $\rho_0$ would be enough to argue as in the previous proof in order to infer the continuity of $\widetilde{U}$.
\end{obs}

\begin{propos} \label{prop:existence of subsolutions}
Let $\rho_0\in C^1_p(Q;\R^+)$ be a given density function as in Proposition \ref{prop:criterion} and let $T$ be any given positive time.
There exist triples $(\overline{m},\overline{U},\overline{q})$ solving
\eqref{eq:linearized system} distributionally on $\R^n\times \R$ enjoying the
following properties:
\begin{align}
&(\overline{m},\overline{U},\overline{q}) \text{ is continuous in } \R^n \times (\R \backslash  \{0\})
\text{ and } \overline{m} \in C(\R ;H_w(\R^n)), \label{eq:thesis 1 subsolution}\\
&\overline{U}(\cdot, t)= \widetilde{U}(\cdot) \text{ for } t=-T,T  \label{eq:thesis 2.2 subsolution}\\
\intertext{and} 
&\overline{q}(x)=p(\rho_0(x))+\frac{\widetilde{\chi}(t)}{n} \text{ for all } (x,t)\in \R^n\times \R, \label{eq:thesis 3 subsolution}\\
& e(\rho_0(x),\overline{m}(x,t),\overline{U}(x,t))< \frac{\widetilde{\chi}(t)}{n}\text{ for all } (x,t)\in \R^n\times ([-T,0[\cap]0,T]). \label{eq:thesis 4 subsolution} 
\end{align}
Moreover
\begin{equation} \label{eq:thesis 5 subsolution}
\abs{\overline{m}(x,0)}^2=\rho_0(x)\chi(0) \text{ a.e. in } \R^n.
\end{equation}
\end{propos}

\begin{dimo}[Proposition \ref{prop:existence of subsolutions}]
We first choose $\overline{q}:=\tilde{q}$ given by Proposition \ref{prop:starting subsolution}. This choice
already yields \eqref{eq:thesis 3 subsolution}.

Now, in analogy with Definition \ref{def:subsolutions} we consider
the space $X_0$ defined as the set of continuous vector fields $m:\R^n\times
]-T,T[ \rightarrow\R^n$ in $C^0(]-T,T[ ;C^0_p (Q) )$ to which
there exists a continuous space-periodic matrix field $U: \R^n \times ]-T,T[\rightarrow
\Sym^n_0$ such that
\begin{align} \label{eq:linearized equations}
    \div_x m = 0,\notag \\
    \dt m +\div_x U+\nabla_x \overline{q}=0,
\end{align}

\begin{align} \label{eq:support}
&\text{supp} (m) \subset Q \times [-T/2,T/2[ \\
&U(\cdot, t)= \widetilde{U}(\cdot) \text{ for all } t\in [-T, T[ \backslash [-T/2, T/2] 
\end{align}
and
\begin{equation} \label{eq:subsolutions constraint}
e(\rho_0(x),m(x,t),U(x,t))<\frac{\widetilde{\chi}(t)}{n}\text{
for all } (x,t)\in \R^n\times ]-T, T[.
\end{equation}
As in Section 4.1, $X_0$ consists of functions $m: ]-T,
T[\rightarrow H$ taking values in a bounded set $B\subset
H$. On $B$ the weak topology of $L^2$ is metrizable, and
correspondingly we find a metric $d$ on $C(]-T, T[; B)$ inducing the
topology of $C(]-T, T[; H_w(Q))$.

Next we note that with minor modifications the proof of Lemma \ref{lemma:improvement lemma} leads to the following claim:

\textbf{Claim:}  Let $Q_0\subset Q$ be given. Let $m\in X_0$ with associated matrix field $U$ and let
$\alpha >0$ such that
$$\int_{Q_0} \left[ \abs{m(x,0)}^2- (\rho_0(x)\widetilde{\chi}(0))\right] dx < -\alpha$$
Then, for any $\delta >0$ there exists a sequence $m_k \in X_0$ with associated smooth matrix field $U_k$ such that
\begin{equation*} 
 \text{supp} (m_k-m, U_k- U) \subset Q_0\times [-\delta, \delta ],
\end{equation*}
\begin{equation*} \label{eq:claim 2}
 m_k \stackrel{d}{\rightarrow} m,
\end{equation*}
and 
\begin{equation*} 
 \liminf_{k\rightarrow\infty} \int_{Q_0} \abs{m_k(x,0)}^2 \geq \int_{Q_0} \abs{m(x,0)}^2  dx +\beta \alpha^2.
\end{equation*}

Fix an exhausting sequence of bounded open subsets
$Q_k\subset Q_{k+1}\subset Q$, each compactly contained
in $\Omega$, and such that $\abs{Q_{k+1}\backslash
Q_k}\leq 2^{-k}$. Let also $\gamma_\varepsilon$ be a standard
mollifying kernel in $\R^n$ (the unusual notation $\gamma_\varepsilon$ for the standard mollifying kernel
is aimed at avoiding confusion between it and the density function).
Using the claim above we construct
inductively a sequence of momentum vector fields $m_k\in X_0$, associated
matrix fields $U_k$ and a sequence of
numbers $\eta_k<2^{-k}$ as follows.

First of all let $m_1 \equiv 0$, $U_1 (x,t)=\widetilde{U}(x)$ for all $(x,t)\in \R^{n+1}$ 
and having obtained $(m_1, U_1)$, ..., $(m_k, U_k)$,
$\eta_1, ..., \eta_{k -1}$ we choose $\eta_k<2^{-k}$ in such a way
that
\begin{equation} \label{eq:recursive step}
\norm{m_k-m_k\ast \gamma_{\eta_k}}_{L^2} < 2^{-k}.
\end{equation}
Then, we set 
$$\alpha_k= -\int_{Q_k} [\abs{m_k (x,0)}^2 - \rho_0 (x)\widetilde{\chi}(0))] dx.$$

Note that \eqref{eq:subsolutions constraint} ensures $\alpha_k>0$.
Then, we apply the claim with $Q_k$, $\alpha= \alpha_k$ and $\delta= 2^{-k} T$ to obtain $m_{k+1}\in X_0$ 
and associated smooth matrix field $U_{k+1}$ such that

\begin{equation} \label{eq:claim 1}
 \text{supp} (m_{k+1}-m_k, U_{k+1}- U_k) \subset Q_k \times [-2^{-k} T, 2^{-k} T],
\end{equation}

\begin{equation} \label{eq:claim 22}
 d(m_{k+1}, m_k) < 2^{-k},
\end{equation}

\begin{equation} \label{eq:claim 3}
 \int_{Q_k} \abs{m_{k+1}(x,0)}^2 dx \geq \int_{Q_k} \abs{m_k(x,0)}^2 dx +\beta \alpha_k^2.
\end{equation}
Since $d$ induces the topology of $C(]-T, T[; H_w(\Omega))$
we can also require that
\begin{equation} \label{eq:claim 4}
 \norm{(m_k-m_{k+1}) \ast \gamma_{\eta_j}}_{L^2(\Omega)} < 2^{-k} \text{ for all } j\leq k \text{ for } t=0.
\end{equation}
From \eqref{eq:claim 2} we infer the existence of a function $\overline{m} \in C(]-T,T[, H_w(\Omega))$ such that
$$m_k\stackrel{d}{\rightarrow} \overline{m}. $$
Besides, \eqref{eq:claim 1} implies that for any compact subset $S$ of $Q\times ]-T,0[ \cup ]0,T[$ there
exists $k_0$ such that $(m_k, U_k)\arrowvert_S =(m_{k_0}, U_{k_0})\arrowvert_S$ for all $k>k_0$. Hence $(m_k, U_k)$
converges in $C^0_{\text{loc}} (Q\times ]-T,0[ \cup ]0,T[)$ to a continuous pair $(\overline{m}, \overline{U})$ solving 
equations \eqref{eq:linearized equations} in $\R^n\times  ]-T,0[ \cup ]0,T[$ and such that \eqref{eq:thesis 1 subsolution}-\eqref{eq:thesis 4 subsolution}
hold. 
In order to conclude, we show that also \eqref{eq:thesis 5 subsolution} holds for $\overline{m}$.

As first, we observe that \eqref{eq:claim 3} yields
$$\alpha_{k+1}\leq \alpha_k -\beta \alpha_k^2+ \abs{Q_{k+1} \backslash Q_{k}} \leq \alpha_k -\beta \alpha_k^2+ 2^{-k},$$
from which we deduce that
$$\alpha_k \rightarrow 0 \text{ as } k\rightarrow \infty.$$
This, together with the following inequality
$$ 0 \geq \int_{Q} \left[\abs{m_k(x,0)}^2 - \rho_0(x)\chi(0) \right] dx \geq -(\alpha_k +C \abs{Q \backslash Q_k}) \geq -(\alpha_k+ C 2^{-k}),$$
implies that
\begin{equation} \label{eq:lim equality}
 \lim_{k\uparrow \infty} \int_{\Omega} \left[\abs{m_k(x,0)}^2 - \rho_0(x)\chi(0) \right] dx =0.
\end{equation}
On the other hand, owing to \eqref{eq:recursive step} and \eqref{eq:claim 4}, we can write for $t=0$ and for every $k$
\begin{align} \label{eq:convolution estimates}
&\norm{m_k - \overline{m}}_{L^2} \notag \\
&\leq \norm{m_k-m_k\ast \gamma_{\eta_k}}_{L^2}+\norm{m_k\ast \gamma_{\eta_k}- \overline{m}\ast \gamma_{\eta_k}}_{L^2}+\norm{\overline{m}\ast \gamma_{\eta_k}- \overline{m}}_{L^2} \notag \\ 
&\leq 2^{-k}+ \sum_{j=0}^{\infty} \norm{m_{k+j} \ast \gamma_{\eta_k}- m_{k+j+1}\ast \gamma_{\eta_k}}_{L^2}+2^{-k} \notag \\
& \leq 2^{-(k-2)}.
\end{align}

Finally, \eqref{eq:convolution estimates} implies that $m_k(\cdot, 0) \rightarrow \overline{m}(\cdot, 0)$ strongly in $H(Q)$ as $k\rightarrow \infty$, which together with \eqref{eq:lim equality}
gives
$$\abs{\overline{m}(x,0)}^2=\rho_0(x)\chi(0) \text{ for almost every } x\in\R^n.$$
\end{dimo}

\section{Proof of the main Theorems}

\begin{dimo} [Proof of Theorem \ref{th:semistationary non-uniqueness}]
Let $T$ be any finite positive time and $\rho_0\in C^1_p(Q)$ be a given density function.
Let also $(\overline{m},\overline{U}, \overline{q})$ be as in Proposition \ref{prop:existence of subsolutions}.
Then, define  $\chi(t):=\widetilde{\chi}(t)$, $q_0(x):=\overline{q}(x)$,
\begin{equation}\label{eq:initial datum for m}
m_0(x,t)=\left\{\begin{array}{ll}
\overline{m}(x,t)  \text{ for } t\in [0,T]\\
\overline{m}(x,t-2T)  \text{ for } t\in [T,2T],
\end{array}\right.
\end{equation}
\begin{equation}\label{eq:initial datum for U}
U_0(x,t)=\left\{\begin{array}{ll}
\overline{U}(x,t)  \text{ for } t\in [0,T]\\
\overline{U}(x,t-2T)  \text{ for } t\in [T,2T].
\end{array}\right.
\end{equation}

For this choices, the quadruple $(m_0, U_0, q_0, \chi)$ satisfies the assumptions of Proposition \ref{prop:criterion}.
Therefore, there exist infinitely many solutions $m\in C([0,2T], H_w(Q))$ of \eqref{eq:semistationary system} in $\R^n\times [0,2T[$ 
with density $\rho_0$, such that
$$m(x,0)=\overline{m}(x,0)=m(x,2T)\text{ for a.e. } x\in\Omega$$
and
\begin{equation} \label{eq:proof 1}
 \abs{m(\cdot,t)}^2= \rho_0(\cdot)\chi(0) \text{ for almost every } (x,t)\in \R^n \times ]0,2T[.
\end{equation}
Since $\abs{m_0(\cdot, 0)}^2=\rho_0(\cdot)\chi(0)$ a.e. in $\R^n$ as well, it is enough to define 
$m^0(x)=m_0(x,0)$ to satisfy also \eqref{eq:modulus constraint b} and hence conclude the proof.
\end{dimo}

\begin{dimo}[Proof of Theorem \ref{th:function chi}]
Under the assumptions of Theorem \ref{th:semistationary non-uniqueness}, we have proven the existence 
of a bounded initial momentum $m^0$ allowing for infinitely many solutions $m\in C([0,T]; H_w(Q))$ of \eqref{eq:semistationary system}
on $\R^n \times [0,T[$ with density $\rho_0$. Moreover, the proof (see Proof of Proposition \ref{prop:starting subsolution}) showed that 
for any smooth function $\chi:\R\rightarrow \R^+$ with $\chi> n\widetilde{\lambda}>0$ the following holds
\begin{align} 
 &\abs{m(x,t)}^2=\rho_0(x) \chi(t) \quad \text{a.e. in }\R^n \times [0,T[, \label{eq:modulus constraint1}\\
&\abs{m^0(x)}^2=\rho_0(x) \chi(0) \quad \text{a.e. in }\R^n. \label{eq:modulus constraint2}
\end{align}
Now, we claim that there exist constants $C_1,C_2>0$ such that choosing the function $\chi(t)>n \widetilde{\lambda}$ on $[0,T[$
among solutions of the following differential inequality
\begin{equation} \label{eq:differential inequality}
 \chi'(t)\leq -C_1 \chi^{1/2}(t) -C_2 \chi^{3/2}(t),
\end{equation}
then the weak solutions $(\rho_0, m)$ of \eqref{eq:semistationary system} obtained in Theorem \ref{th:semistationary non-uniqueness}
will also satisfy the admissibility condition \eqref{eq:admissibility condition} on $\R^n\times [0,T[$.
Of course, there is an issue of compatibility between the differential inequality \eqref{eq:differential inequality} and 
the condition $\chi> n \widetilde{\lambda}$: this motivates the existence of a time $\overline{T}>0$ defining the maximal time-interval in which
the admissibility condition indeed holds.

Let $T$ be any finite positive time. As first, we aim to prove the claim.
Since $m\in C([0,T];H_w(Q))$ is divergence-free and fulfills \eqref{eq:modulus constraint1}-\eqref{eq:modulus constraint2} 
and $\rho_0$ is time-independent, 
\eqref{eq:admissibility condition} reduces to the following inequality
\begin{equation} \label{eq:admissibility condition reduced}
 \frac{1}{2} \chi'(t) + m\cdot \nabla \Bigg( \varepsilon (\rho_0(x)) 
+\frac{p(\rho_0(x))}{\rho_0(x)}\Bigg)+\frac{\chi(t)}{2} m \cdot \nabla\left(\frac{1}{\rho_0(x)}\right) \leq 0,
\end{equation}
intended in the sense of (space-periodic) distributions on $\R^n \times [0,T]$.
As $\rho_0 \in C^1_p(Q)$, there exists a constant $c_0^2$ with $\rho_0\leq c_0^2$ on $\R^n$, 
whence (see \eqref{eq:modulus constraint1}-\eqref{eq:modulus constraint2} )
\begin{equation}\label{eq:bound for modulus of m}
 \abs{m(x,t)}\leq c_0 \sqrt{\chi(t)} \text{ a.e. on } \R^n \times [0,T[.
\end{equation}
Similarly we can find constants $c_1, c_2>0$ with
\begin{align}
 &\abs{\nabla \Bigg( \varepsilon (\rho_0(x)) 
+\frac{p(\rho_0(x))}{\rho_0(x)}\Bigg)}\leq c_1 \text{  a.e. in }\R^n \label{eq:first constant} \\
&\abs{\nabla \left(\frac{1}{\rho_0(x)}\right)}\leq c_2 \text{  a.e. in }\R^n. \label{eq:second constant}
\end{align}
As a conseguence of \eqref{eq:bound for modulus of m}-\eqref{eq:second constant}, \eqref{eq:admissibility condition reduced}
holds as soon as $\chi$ satisfies
\begin{equation*}
  \chi'(t)\leq -2 c_1 c_0 \chi^{1/2}(t) -c_2 c_0 \chi^{3/2}(t) \text{ on } [0,T[.
\end{equation*}
Therefore, by choosing $C_1:=2 c_1 c_0$ and $C_2:=c_2 c_0$ we can conclude the proof of the claim.

Now, it remains to show the existence of a function $\chi$ as in the claim, i.e. that both the differential inequality 
\eqref{eq:differential inequality} and the condition $\chi> n\widetilde{\lambda}$ can hold true on some suitable time-interval.
To this aim, we can consider the equality in \eqref{eq:differential inequality}, couple it with
the initial condition $\chi(0)= \chi_0$ for some constant $\chi_0 >n\widetilde{\lambda}$ and then solve the resulting Cauchy problem.
For the obtained solution $\chi$, there exists a positive time $\overline{T}$ such that 
$\chi(t) > n \widetilde{\lambda}$ on $[0,\overline{T}[$. 

Finally, applying the claim on the time-interval $[0,\overline{T}[$ we conclude that the admissibility condition holds on
$\R^n \times [0,\overline{T}[$ as desired.
\end{dimo}

\begin{dimo} [Proof of Theorem \ref{th:isentropic non-uniqueness}]
 The proof of Theorem \ref{th:isentropic non-uniqueness} strongly relies on Theorems
\ref{th:semistationary non-uniqueness}-\ref{th:function chi}.
Given a continuously differentiable initial density $\rho^0$ we apply Theorems 
\ref{th:semistationary non-uniqueness}-\ref{th:function chi} for $\rho_0(x):=\rho^0(x)$ thus obtaining a positive time $\overline{T}$ and
a bounded initial momentum $m^0$ allowing for infinitely many solutions $m\in C([0,T]; H_w(Q))$ of \eqref{eq:semistationary system}
on $\R^n \times [0,\overline{T}[$ with density $\rho^0$ and such that the following holds
\begin{align} 
 &\abs{m(x,t)}^2=\rho_0(x) \chi(t) \quad \text{a.e. in }\R^n \times [0,\overline{T}[, \label{eq:modulus constraint11}\\
&\abs{m^0(x)}^2=\rho_0(x) \chi(0) \quad \text{a.e. in }\R^n, \label{eq:modulus constraint22}
\end{align}
for a suitable smooth function $\chi:[0,\overline{T}]\rightarrow \R^+$.
Now, define $\rho(x,t)=\rho_0(x)\mathbf{1}_{[0, \overline{T}[}(t)$.
This shows that \eqref{eq:weak2} holds.
To prove \eqref{eq:weak1} observe that $\rho$ is independent of $t$ and $m$ is weakly divergence-free for almost every $0<t<\overline{T}$.
Therefore, the pair $(\rho, m)$ is a weak solution of \eqref{eq:isentropic Euler system} with initial data $(\rho^0, m^0)$. 
Finally, we can also prove \eqref{eq:admissibility condition}: each solution obtained is also admissible.
Indeed, for $\rho(x,t)=\rho_0(x)\mathbf{1}_{[0, \overline{T}[}(t)$, \eqref{eq:admissibility condition} is ensured by Theorem \ref{th:function chi}.
\end{dimo}

\section*{Acknowledgements}
The author would like to thank Professor Camillo De Lellis for his support and for many helpful conversations.

\end{document}